\documentclass[12pt,oneside]{article}
\usepackage{amsmath,amssymb,amscd,latexsym,amsthm}
\usepackage[unicode]{hyperref}
\usepackage{hypbmsec}
\ifx\pdfoutput\undefined\else
\advance\topmargin by -15mm
\hypersetup{
colorlinks=true,
linkcolor=blue,
citecolor=blue,
urlcolor=blue,
filecolor=blue,
bookmarksnumbered=true,
pdfstartview=FitH,
pdfhighlight=/N
}
\fi
\makeatletter
\def\ps@myheadings{%
    \let\@oddfoot\@empty\let\@evenfoot\@empty
    \def\@evenhead{\small\thepage\hfil{\itshape\leftmark}\hfil\phantom{\thepage}}%
    \def\@oddhead{\small\phantom{\thepage}\hfil{\itshape\rightmark}\hfil\thepage}%
    \let\@mkboth\@gobbletwo
    \let\sectionmark\@gobble
    \let\subsectionmark\@gobble
    }
\makeatother
\numberwithin{equation}{section}
\renewcommand{\d}{\mathrm d}
\newcommand{\ol}{\overline}
\newcommand{\wt}{\widetilde}

\DeclareMathOperator{\Imag}{Im}
\renewcommand{\Im}{\Imag}
\begin{document}
\newtheorem{theorem}{Theorem}
\theoremstyle{remark}
\newtheorem{example}{Example}
\newtheorem{conj}{Conjecture}
\newtheorem*{remark}{Remark}

\hypersetup{pdfauthor={Wadim Zudilin},%
pdftitle={Ramanujan-type formulae for \$1/\000\134pi\$: A second wind?}}

\title{Ramanujan-type formulae for $1/\pi$: A~second~wind?%
\footnote{The work was supported by a fellowship
of the Max Planck Institute for Mathematics (Bonn) and supported
in part by the INTAS foundation, grant no.~03-51-5070.\newline
2000 {\it AMS subject classification\/}: 11F11, 11Y60, 33C20 (primary),
05A19, 11B65, 11J82, 11M06, 14H52, 14J32, 33C75, 33F10, 34M50, 40G99, 65B10, 65Q05 (secondary).}}

\author{Wadim Zudilin}

\date{14 May 2008}

\pagestyle{myheadings}
\markboth{Wadim Zudilin}{Ramanujan-type formulae for $1/\pi$\/\textup: A~second~wind\/\textup?}

\maketitle

\begin{abstract}
In 1914 S.~Ramanujan recorded a list of 17 series for $1/\pi$.
We survey the methods of proofs of Ramanujan's formulae and indicate
recently discovered generalizations, some of which are not yet proven.
\end{abstract}

\bigskip

The twentieth century was full of mathematical discoveries.
Here we expose two significant contributions from that time,
in reverse chronological order.
At first glance, the stories might be thought of a different nature.
But we will try to convince the reader that they have much in common.

\section({Ramanujan and Ap\'ery: \$1/\000\134pi\$ and \$\000\134zeta(3)\$})%
{Ramanujan and Ap\'ery: $1/\pi$ and $\zeta(3)$}
\label{s1}

In 1978 R.~Ap\'ery showed the irrationality of
$\zeta(3)$ (see \cite{Ap} and~\cite{Po}).
His rational approximations to the number in question
(known nowadays as \emph{Ap\'ery's constant\/})
have the form $v_n/u_n\in\mathbb Q$ for $n=0,1,2,\dots$,
where the denominators $\{u_n\}=\{u_n\}_{n=0,1,\dots}$
and numerators $\{v_n\}=\{v_n\}_{n=0,1,\dots}$ satisfy
the same polynomial recurrence
\begin{equation}
(n+1)^3u_{n+1}-(2n+1)(17n^2+17n+5)u_n+n^3u_{n-1}=0
\label{e00}
\end{equation}
with the initial data
\begin{equation*}
u_0=1, \quad u_1=5, \qquad v_0=0, \quad v_1=6.
\end{equation*}
Then
\begin{equation*}
\lim_{n\to\infty}\frac{v_n}{u_n}=\zeta(3)
\end{equation*}
and, surprisingly, the denominators $\{u_n\}$ are integers:
\begin{equation}
u_n=\sum_{k=0}^n\binom nk^2\binom{n+k}k^2\in\mathbb Z,
\qquad n=0,1,2,\dots,
\label{e01}
\end{equation}
while the numerators $\{v_n\}$ are `close' to being integers.

In 1914 S.~Ramanujan~\cite{Ra}, \cite{Be} recorded a list of 17 series
for $1/\pi$, from which we indicate the simplest one
\begin{equation}
\sum_{n=0}^\infty\frac{(\frac12)_n^3}{n!^3}
(4n+1)\cdot(-1)^n
=\frac2\pi
\label{e02}
\end{equation}
and also two quite impressive examples
\begin{align}
\sum_{n=0}^\infty\frac{(\frac14)_n(\frac12)_n(\frac34)_n}{n!^3}
(21460n+1123)\cdot\frac{(-1)^n}{882^{2n+1}}
&=\frac4\pi\,,
\label{e03}
\\
\sum_{n=0}^\infty\frac{(\frac14)_n(\frac12)_n(\frac34)_n}{n!^3}
(26390n+1103)\cdot\frac1{99^{4n+2}}
&=\frac1{2\pi\sqrt2}
\label{e04}
\end{align}
which produce rapidly converging (rational) approximations to~$\pi$. Here
$$
(a)_n=\frac{\Gamma(a+n)}{\Gamma(a)}=\begin{cases}
a(a+1)\dotsb(a+n-1) &\text{for $n\ge1$}, \\
1 &\text{for $n=0$},
\end{cases}
$$
denotes the Pochhammer symbol (the rising factorial). The Pochhammer
products occurring in all formulae of this type may be written
in terms of binomial coefficients:
\begin{gather*}
\frac{(\frac12)_n^3}{n!^3}
=2^{-6n}\binom{2n}n^3,
\qquad
\frac{(\frac13)_n(\frac12)_n(\frac23)_n}{n!^3}
=2^{-2n}3^{-3n}\binom{2n}n\frac{(3n)!}{n!^3},
\\
\frac{(\frac14)_n(\frac12)_n(\frac34)_n}{n!^3}
=2^{-8n}\frac{(4n)!}{n!^4},
\qquad
\frac{(\frac16)_n(\frac12)_n(\frac56)_n}{n!^3}
=12^{-3n}\frac{(6n)!}{n!^3(3n)!}.
\end{gather*}
Ramanujan's original list was subsequently extended to
several other series which we plan to
touch on later in the paper. For the moment
we give two more celebrated examples:
\begin{align}
\sum_{n=0}^\infty\frac{(\frac13)_n(\frac12)_n(\frac23)_n}{n!^3}
(14151n+827)\cdot\frac{(-1)^n}{500^{2n+1}}
&=\frac{3\sqrt3}\pi\,,
\label{e05}
\\
\sum_{n=0}^\infty\frac{(\frac16)_n(\frac12)_n(\frac56)_n}{n!^3}
(545140134n+13591409)\cdot\frac{(-1)^n}{53360^{3n+2}}
&=\frac3{2\pi\sqrt{10005}}\,.
\label{e06}
\end{align}
Formula~\eqref{e05} is proven by H.\,H.~Chan, W.-C.~Liaw and V.~Tan
\cite{CLT} and \eqref{e06} is the Chudnovskys' famous formula~\cite{Ch} which enabled
them to hold the record for the calculation of~$\pi$ in 1989--94.
On the left-hand side of each formula~\eqref{e02}--\eqref{e06}
we have linear combinations of a (generalized) hypergeometric series
\begin{equation}
{}_mF_{m-1}\biggl(\begin{matrix}
a_1, & a_2, & \dots, & a_m \\
& b_2, & \dots, & b_m \end{matrix} \biggm|z\biggr)
=\sum_{n=0}^\infty
\frac{(a_1)_n(a_2)_n\cdots(a_m)_n}
{(b_2)_n\cdots(b_m)_n}\,
\frac{z^n}{n!}
\label{e07}
\end{equation}
and its derivative at a point close to the origin.
The rapid convergence of the series in \eqref{e03}--\eqref{e06}
may be used for proving the quantitative irrationality
of the numbers $\pi\sqrt d$ with $d\in\mathbb N$ (see \cite{Z1}
for details).

In both Ramanujan's and Ap\'ery's cases, there were just
hints on how the things might be proven. Rigorous proofs appeared
somewhat later. We will not discuss proofs of Ap\'ery's
theorem and its further generalizations (see \cite{Fi}
for a review of the subject), just concentrating
on the things around the remarkable Ramanujan-type series.
But we will see that both Ramanujan's and Ap\'ery's discoveries
have several common grounds.

\section{Elliptic proof of Ramanujan's formulae}
\label{s2}

Although Ramanujan did not indicate how he arrived at his series,
he hinted that these series belong to what is now known as
`the theories of elliptic functions to alternative bases'.
The first rigorous mathematical proofs of Ramanujan's series and
their generalizations were given by the Borweins~\cite{BB} and
Chudnovskys~\cite{Ch}. Let us sketch, following~\cite{Ch},
the basic ideas of those very first proofs.

One starts with an elliptic curve $y^2=4x^3-g_2x-g_3$
over $\ol{\mathbb Q}$ with fundamental periods $\omega_1,\omega_2$
(where $\Im(\omega_2/\omega_1)>0$) and corresponding
quasi-periods $\eta_1,\eta_2$. Besides the Legendre relation
\begin{equation*}
\eta_1\omega_2-\eta_2\omega_1=2\pi i,
\end{equation*}
the following \emph{linear\/}
relations between $\omega_1,\omega_2,\eta_1,\eta_2$
over~$\ol{\mathbb Q}$ are available in the complex multiplication case,
i.e., when $\tau=\omega_2/\omega_1\in\mathbb Q[\sqrt{-d}]$
for some $d\in\mathbb N$:
\begin{equation}
\omega_2-\tau\omega_1=0, \qquad
A\tau\eta_2-C\eta_1+(2A\tau+B)\alpha\omega_1=0,
\label{e08}
\end{equation}
where the integers $A$, $B$ and $C$ come from the equation
$A\tau^2+B\tau+C=0$ defining the quadratic number~$\tau$
and $\alpha\in\mathbb Q(\tau,g_2,g_3)\subset\ol{\mathbb Q}$.
Equations~\eqref{e08} allow one to express $\omega_2,\eta_2$
by means of $\omega_1,\eta_1$ only. Substituting these
expressions into~\eqref{e08}, and using the hypergeometric
formulae for $\omega_1,\eta_1$ and also for
$\omega_1^2,\omega_1\eta_1$ (which follow from Clausen's identity)
one finally arrives at a formula of Ramanujan type.
An important (and complicated) problem in the proof is
computing the algebraic number
$$
\alpha=\frac{8\pi^2}{81\omega_1^2}\biggl(E_2(\tau)-\frac3{\pi\Im\tau}\biggr),
\qquad\text{where}\quad
E_2(\tau)=1-24\sum_{n=1}^\infty e^{2\pi in\tau}\sum_{d\mid n}d.
$$
Note that $\alpha$ viewed as a function of~$\tau$
is a non-holomorphic modular form of weight~2. Although the Chudnovskys
attribute the knowledge of the fact that $\alpha(\tau)$
takes values in the Hilbert class field $\mathbb Q(\tau,j(\tau))$
of~$\mathbb Q(\tau)$ to Kronecker (in Weil's
presentation~\cite{We}), we would refer the reader to the work~\cite{BC}
by B.\,C.~Berndt and H.\,H.~Chan.

\section{Modular proof of Ramanujan's formulae}
\label{s3}

An understanding of the complication of the above proof came in 2002
with T.~Sato's discovery of the formula
\begin{equation}
\sum_{n=0}^\infty u_n\cdot(20n+10-3\sqrt5)
\biggl(\frac{\sqrt5-1}2\biggr)^{12n}
=\frac{20\sqrt3+9\sqrt{15}}{6\pi}
\label{e09}
\end{equation}
of Ramanujan type, involving Ap\'ery's numbers~\eqref{e01}.
The modular argument was essentially simplified by H.\,H.~Chan
with his collaborators and later by Y.~Yang to produce a lot of new
identities like~\eqref{e09} based on a not necessarily hypergeometric
series $F(z)=\sum_{n=0}^\infty u_nz^n$. Examples are
\begin{equation}
\sum_{n=0}^\infty\sum_{k=0}^n\binom nk^2\binom{2k}k\binom{2n-2k}{n-k}
\cdot(5n+1)\frac1{64^n}
=\frac8{\pi\sqrt3}
\label{e10}
\end{equation}
due to H.\,H.~Chan, S.\,H.~Chan and Z.-G.~Liu~\cite{CCL};
\begin{equation}
\sum_{n=0}^\infty\sum_{k=0}^{[n/3]}(-1)^{n-k}3^{n-3k}
\frac{(3k)!}{k!^3}\binom n{3k}\binom{n+k}k
\cdot(4n+1)\frac1{81^n}
=\frac{3\sqrt3}{2\pi}
\label{e11}
\end{equation}
due to H.\,H.~Chan and H.~Verrill (2005);
\begin{equation}
\sum_{n=0}^\infty\sum_{k=0}^n\binom nk^4
\cdot(4n+1)\frac1{36^n}
=\frac{18}{\pi\sqrt{15}}
\label{e12}
\end{equation}
due to Y.~Yang (2005).

It should be mentioned that Picard--Fuchs differential
equations (of order~3) satisfied by the series $F(z)$ always have very
nice arithmetic properties~\cite{Ya}. Therefore, it is not surprising
that $F(z)$ admits a modular parametri\-zation: $f(\tau)=F(z(\tau))$
is a modular form of weight~2 for a modular (uniformizing)
substitution $z=z(\tau)$.

Let us follow Yang's argument to show the basic ideas of
the new proof in the example of~\eqref{e09}. Our choice is
$$
z(\tau)=\biggl(\frac{\eta(\tau)\eta(6\tau)}
{\eta(2\tau)\eta(3\tau)}\biggr)^{12},
\qquad
f(\tau)=\frac{\eta(2\tau)^7\eta(3\tau)^7}
{\eta(\tau)^5\eta(6\tau)^5},
$$
which are modular forms of level~6; the expressions were
obtained by F.~Beukers in his proof of Ap\'ery's theorem
using modular forms~\cite{Beu}. Here
$$
\eta(\tau)=e^{\pi i\tau/12}\prod_{n=1}^\infty(1-e^{2\pi in\tau})
$$
is the Dedekind eta-function. The function
$g(\tau)=(2\pi i)^{-1}f'(\tau)/f(\tau)$ satisfies
the functional equation
$$
g(\gamma\tau)=\frac{c(c\tau+d)}{\pi i}+(c\tau+d)^2g(\tau)
\qquad\text{for}\quad
\gamma=\begin{pmatrix} a &\ b \\ c &\ d\end{pmatrix}
\in\Gamma_0(6)+w_6,
$$
where $w_6$ denotes the Atkin--Lehner involution.
Taking
$$
\gamma=\frac1{\sqrt6}\begin{pmatrix} 0 &\ -1 \\ 6 &\ 0\end{pmatrix},
\qquad
\tau=\tau_0=\frac i{\sqrt{30}}
$$
we obtain
\begin{equation}
g(\tau_0)+5g(5\tau_0)=\frac{\sqrt{30}}\pi.
\label{e13}
\end{equation}
On the other hand, $h(\tau)=g(\tau)-5g(5\tau)$ is a modular
form of weight~2 and level~30
(on $\Gamma_0(30)+\langle w_5,w_6\rangle$).
This implies that $h(\tau)/f(\tau)$ is an algebraic
function of~$z(\tau)$ and after explicit evaluations at
$\tau=\tau_0$ we arrive at
\begin{equation}
g(\tau_0)-5g(5\tau_0)=h(\tau_0)
=\frac{900\sqrt2-402\sqrt{10}}5f(\tau_0)
=(900\sqrt2-402\sqrt{10})f(5\tau_0).
\label{e14}
\end{equation}
Combining \eqref{e13} and \eqref{e14} we deduce that
$$
\frac{\sqrt{30}}\pi
=(900\sqrt2-402\sqrt{10})f(5\tau_0)+10g(5\tau_0),
$$
and it only remains to use the expansion
$$
g(5\tau_0)
=z\frac{\d f/\d z}f\cdot
\frac1{2\pi i}\frac{z'(\tau)}{z(\tau)}\bigg|_{\tau=5\tau_0}
=(108\sqrt2-48\sqrt{10})\sum_{n=0}^\infty nu_n\cdot z(5\tau_0)^n
$$
(since $z'(\tau)/(2\pi i)$ and $f(\tau)$ are modular forms
of weight~$2$ on~$\Gamma_0(6)+w_6$, the function
$$
\frac1f\cdot
\frac1{2\pi i}\frac{z'}z
$$
is an algebraic function of~$z$)
and the evaluation
$$
z(5\tau_0)=z(\tau_0)=161-72\sqrt5=\biggl(\frac{\sqrt5-1}2\biggr)^{12}.
$$

As pointed out to us by H.\,H.~Chan, the main difficulty one meets
in the above proof is to prove the algebraicity evaluations rigorously
(cf.~\cite{BC}).

\section{Creative telescoping}
\label{s4}

There is yet another method of proof, but applicable only to a small number of
Rama\-nu\-jan-type series. It is based on the algorithm of creative
telescoping, due to Gosper--Zeilberger. Note that
an essential part of the first proof of Ap\'ery's theorem~\cite{Po}, namely,
the proof of the recurrence~\eqref{e00}, was given by D.~Zagier
also using a telescoping argument. D.~Zeilberger
(and his automatic collaborator S.\,B.~Ekhad) could prove the simplest
Ramanujan's identity~\eqref{e02} in the following way~\cite{EZ}.
One verifies the (terminating) identity
\begin{equation}
\sum_{n=0}^\infty\frac{(1/2)_n^2(-k)_n}{n!^2(3/2+k)_n}
(4n+1)(-1)^n
=\frac{\Gamma(3/2+k)}{\Gamma(3/2)\Gamma(1+k)}
\label{e15}
\end{equation}
for all \emph{non-negative\/} integers $k$. To do this,
divide both sides of~\eqref{e15} by the right-hand side
and denote the summand on the left by $F(n,k)$:
$$
F(n,k)=(4n+1)(-1)^n\frac{(1/2)_n^2(-k)_n}{n!^2(3/2+k)_n}\,
\frac{\Gamma(3/2)\Gamma(1+k)}{\Gamma(3/2+k)};
$$
then take
$$
G(n,k)=\frac{(2n+1)^2}{(2n+2k+3)(4n+1)}F(n,k)
$$
with the motive that $F(n,k+1)-F(n,k)=G(n,k)-G(n-1,k)$,
hence $\sum_nF(n,k)$ is a constant, which is seen to be~1
by plugging in $k=0$. Finally, to deduce \eqref{e02} one
takes $k=-1/2$, which is legitimate in view
of Carlson's theorem~\cite[Section~5.3]{Ba}.

\section(Guillera's series for \$1/\000\134pi\000\136{}2\$)%
{Guillera's series for $1/\pi^2$}
\label{s5}

If one wishes to use the latter method of proof for
other Ramanujan-type formulae, ingenuity is required in order to put
the new parameter~$k$ in the right place.
This was done only recently by J.~Guillera \cite{G1,G3},
who used the method to prove some other identities of Ramanujan
(in those cases when $z$ has only 2 and 3 in its prime
decomposition). If the reader doubts the applicability
of the method, then take into account that the purely hypergeometric
origin of the method and its independence from the elliptic and modular stuff
allowed Guillera~\cite{G1,G2,G3} to prove new generalizations
of Ramanujan-type series, namely,
\begin{align}
\sum_{n=0}^\infty\frac{(\frac12)_n^5}{n!^5}
(20n^2+8n+1)\frac{(-1)^n}{2^{2n}}&=\frac8{\pi^2},
\label{e16}
\\
\sum_{n=0}^\infty\frac{(\frac12)_n^5}{n!^5}
(820n^2+180n+13)\frac{(-1)^n}{2^{10n}}&=\frac{128}{\pi^2},
\label{e17}
\\
\sum_{n=0}^\infty\frac{(\frac12)_n^3(\frac14)_n(\frac34)_n}{n!^5}
(120n^2+34n+3)\frac1{2^{4n}}&=\frac{32}{\pi^2},
\label{e18}
\end{align}
and also to find experimentally~\cite{G2} four additional formulae
\begin{align}
\sum_{n=0}^\infty\frac{(\frac12)_n(\frac14)_n(\frac34)_n(\frac16)_n(\frac56)_n}{n!^5}
(1640n^2+278n+15)\frac{(-1)^n}{2^{10n}}&=\frac{256\sqrt3}{3\pi^2},
\label{e19}
\\
\sum_{n=0}^\infty\frac{(\frac12)_n(\frac14)_n(\frac34)_n(\frac13)_n(\frac23)_n}{n!^5}
(252n^2+63n+5)\frac{(-1)^n}{48^n}&=\frac{48}{\pi^2},
\label{e20}
\\
\sum_{n=0}^\infty\frac{(\frac12)_n(\frac13)_n(\frac23)_n(\frac16)_n(\frac56)_n}{n!^5}
(5418n^2+693n+29)\frac{(-1)^n}{80^{3n}}&=\frac{128\sqrt5}{\pi^2},
\label{e21}
\\
\sum_{n=0}^\infty\frac{(\frac12)_n(\frac18)_n(\frac38)_n(\frac58)_n(\frac78)_n}{n!^5}
(1920n^2+304n+15)\frac1{7^{4n}}&=\frac{56\sqrt7}{\pi^2}.
\label{e22}
\end{align}
As Guillera notices, the series in~\eqref{e20}--\eqref{e22}
are closely related to the series
\begin{align*}
\sum_{n=0}^\infty\frac{(\frac12)_n(\frac14)_n(\frac34)_n}{n!^3}
(28n+3)\frac{(-1)^n}{48^n}&=\frac{16}{\pi\sqrt3},
\\
\sum_{n=0}^\infty\frac{(\frac12)_n(\frac16)_n(\frac56)_n}{n!^3}
(5418n+263)\frac{(-1)^n}{80^{3n}}&=\frac{640\sqrt{15}}{3\pi},
\\
\sum_{n=0}^\infty\frac{(\frac12)_n(\frac14)_n(\frac34)_n}{n!^3}
(40n+3)\frac1{7^{4n}}&=\frac{49}{3\pi\sqrt3},
\end{align*}
respectively, proven by the methods in Sections~\ref{s2} and~\ref{s3}.
However, there is no obvious way to deduce
any of formulae \eqref{e16}--\eqref{e22} by modular means;
the problem lies in the fact that the (Zariski closure of
the) projective monodromy group for the corresponding series
$F(z)=\sum_{n=0}^\infty u_nz^n$ is always $O_5(\mathbb R)$
(this is an immediate consequence of a general result
of F.~Beukers and G.~Heckman~\cite{BH}),
which is essentially `richer' than $O_3(\mathbb R)$ for classical
Ramanujan's series.

There exists also the higher-dimensional identity
$$
\sum_{n=0}^\infty\frac{(\frac12)_n^7}{n!^7}
(168n^3+76n^2+14n+1)\frac1{2^{6n}}=\frac{32}{\pi^3},
$$
discovered by B.~Gourevich in 2002 (using an integer
relations algorithm). Guillera also found experimentally
an analogue of Sato's series:
\begin{equation}
\begin{gathered}
\sum_{n=0}^\infty v_n\cdot(36n^2+12n+1)\frac1{2^{10n}}
=\frac{32}{\pi^2},
\\
\text{where}\quad
v_n=\binom{2n}n^2\sum_{k=0}^n\binom{2k}k^2\binom{2n-2k}{n-k}^2.
\end{gathered}
\label{e23}
\end{equation}

\section{Transformations of hypergeometric series}
\label{s6}

As we have seen, Ramanujan's original formulae as well as
Guillera's formulae \eqref{e16}--\eqref{e22} involve classical
hypergeometric series~\eqref{e07}, while series like \eqref{e09}--\eqref{e12}
and~\eqref{e23} are based on double hypergeometric series. A natural way to pass
from one formula to another is by algebraic transformations of the
hypergeometric series involved. For instance, formula~\eqref{e09}
may be deduced from the transformation
\begin{equation*}
\sum_{n=0}^\infty u_nz^n
=\frac1{2+2z-\sqrt{1-34z+z^2}}
\cdot{}_3F_2\biggl(\begin{matrix}
\frac14, & \frac12, & \frac34 \\
& 1, & 1 \end{matrix} \biggm| 256t(z)\biggr)
\end{equation*}
where
\begin{align*}
t(z)
&=\frac z{2(1+14z+z^2)^4}
\bigl(1-36z+199z^2+184z^3+199z^4-36z^5+z^6
\\ &\qquad
+(1+z)(1-z)^2(1-18z+z^2)\textstyle\sqrt{1-34z+z^2}\,\bigr),
\end{align*}
given by Y.~Yang (2005), together with the Ramanujan-type formula for
the ${}_3F_2$-series on the right-hand side specialized at
the point
$$
\biggl(\frac{5+4\sqrt2}{7\sqrt3}\biggr)^4
=256t\biggl(\biggl(\frac{\sqrt5-1}2\biggr)^{12}\biggr).
$$
A similar argument is used by M.\,D.~Rogers in~\cite{Ro}
to deduce some further identities for $1/\pi$ of Ramanujan--Sato type.

Using the quadratic transformation $z\mapsto-4z/(1-z)^2$ of the hypergeometric
series, we were able to produce from~\eqref{e16},~\eqref{e17}
two more series of the latter type~\cite{Z2}:
\begin{align*}
\sum_{n=0}^\infty w_n\frac{(4n)!}{n!^2(2n)!}
(18n^2-10n-3)\frac1{(2^85^2)^n}
&=\frac{10\sqrt{5}}{\pi^2},
\\
\sum_{n=0}^\infty w_n\frac{(4n)!}{n!^2(2n)!}
(1046529n^2+227104n+16032)\frac1{(5^441^2)^n}
&=\frac{5^441\sqrt{41}}{\pi^2},
\end{align*}
where the sequence of integers
$$
w_n=\sum_{k=0}^n\binom{2k}k^3\binom{2n-2k}{n-k}2^{4(n-k)},
\qquad n=0,1,2,\dots,
$$
satisfies the recurrence relation
$$
(n+1)^3w_{n+1}-8(2n+1)(8n^2+8n+5)w_n+4096n^3w_{n-1}=0,
\qquad n=1,2,\dotsc.
$$
In~\cite{Z3} we show that a huge family of formulae for $1/\pi^2$
(as well as for $1/\pi^3$, $1/\pi^4$, etc) can be derived by
taking powers of Ramanujan-type formulae for $1/\pi$.
For instance, the square of the Chudnovskys' formula~\eqref{e06} takes
the monstrous form
\begin{align*}
\sum_{n=0}^\infty w_n\frac{(3n)!}{n!^3}
&(222883324273153467n^2
+16670750677895547n
\\ &\qquad
+415634396862086)\frac{(-1)^n}{640320^{3n+3}}
=\frac1{64\pi^2}.
\end{align*}

\section{Further observations and open problems}
\label{s7}

It is worth mentioning that identities like~\eqref{e15}
are valid for all non-negative \emph{real\/} values of~$k$.
This fact has several other curious implications; for instance,
the series
\begin{equation}
G(k)=\sum_{n=0}^\infty\frac{(1/2+k)_n^5}{(1+k)_n^5}
\bigl(820(n+k)^2+180(n+k)+13\bigr)\frac{(-1)^n}{2^{10n}}
\label{e24}
\end{equation}
has a closed-form evaluation at $k=0$ and $k=1/2$:
$$
G(0)=\frac{128}{\pi^2}
\qquad\text{and}\qquad
G(\tfrac12)=256\zeta(3),
$$
where the first formula follows from~\eqref{e17} while
the second one was given by T.~Amdeberhan and D.~Zeilberger~\cite{AmZe}.
Guillera has conjectured (and proven) evaluations for series
like~\eqref{e24} viewed as functions of the continuous
(complex or real) parameter~$k$.

It seems to be a challenge to develop a modular-like
theory for proving Guillera's identities and finding
a (more or less) general pattern of them.
For the moment, we have only speculations
in this respect on a relationship to mirror symmetry,
namely, to the linear differential equations for the periods
of certain Calabi--Yau threefolds. A standard example here
is the hypergeometric series (cf.~\eqref{e16} and~\eqref{e17})
$$
F(z)={}_5F_4\biggl(\begin{matrix}
\frac12, & \frac12, & \frac12, & \frac12, & \frac12 \\
& 1, & 1, & 1, & 1 \end{matrix} \biggm|2^{10}z\biggr)
=\sum_{n=0}^\infty\binom{2n}n^5z^n,
$$
which satisfies the 5th-order linear differential equation
$$
\bigl(\theta^5-32z(2\theta+1)^5\bigr)Y=0,
\qquad\text{where}\quad \theta=z\frac{\d}{\d z}.
$$
If $G(z)$ is another solution of the latter equation
of the form $F(z)\log z+F_1(z)$ with $F_1(z)\in z\mathbb Q[[z]]$,
then
$$
\wt F(z)=(1-2^{10}z)^{-1/2}\det\begin{pmatrix}
F &\ \ G \\ \theta F &\ \ \theta G \end{pmatrix}^{1/2}
$$
(the sharp normalization factor $(1-2^{10}z)^{-1/2}$ is due to Y.~Yang)
satisfies the 4th-order equation
\begin{equation}
\bigl(\theta^4-16z(128\theta^4+256\theta^3+304\theta^2+176\theta+39)
+2^{20}z^2(\theta+1)^4\bigr)Y=0
\label{e25}
\end{equation}
(entry \#204 in \cite[Table~A]{AESZ}). For a quadratic transformation
of the new function $\wt F(z)$ we have the following explicit formula~\cite{AZ2}:
\begin{equation}
\frac{1+z}{(1-z)^2}\wt F\biggl(\frac{-z}{(1-z)^2}\biggr)
=\sum_{n=0}^\infty\biggl(\sum_{k=0}^n4^{n-k}\binom{2k}k^2\binom{2n-2k}{n-k}\biggr)^2z^n,
\label{e26}
\end{equation}
where the right-hand side is the Hadamard square of the series
\begin{align*}
\frac1{1-16z}{}_2F_1\biggl(\begin{matrix}
\frac12, & \frac12 \\ & 1 \end{matrix}\biggm|\frac{-16z}{1-16z}\biggr)
&=\sum_{n=0}^\infty\binom{2n}n^2\frac{(-1)^nz^n}{(1-16z)^{n+1}}
\\
&=\sum_{n=0}^\infty\biggl(\sum_{k=0}^n4^{n-k}\binom{2k}k^2\binom{2n-2k}{n-k}\biggr)z^n
\end{align*}
which admits a modular uniformization (cf.~Section~\ref{s3}).
It is worth mentioning that \eqref{e25} and the differential equation of order~4
for the right-hand side of~\eqref{e26} are of Calabi--Yau type \cite{AZ}, \cite{CYY},
i.e., they imitate all properties of a differential equation
for the periods of a Calabi--Yau threefold. Are there analogues
of Hilbert class fields for this and similar situations?
Can formulae \eqref{e16}--\eqref{e23} be deduced from
formulae for $1/\pi$ by means of algebraic transformations of
hypergeometric series? There is still some work to do in the subject
originated by Ramanujan's note~\cite{Ra} almost 100~years ago.

\smallskip
\textbf{Acknowledgements.}
It is a pleasure for me to thank Gert Almkvist, Heng Huat Chan,
Jes\'us Guillera, Jonathan Sondow, Yifan Yang, Don Zagier, and
the anonymous referee for helpful comments.

%==================================================

\bigskip
\noindent
\emph{Affiliation}:
Steklov Mathematical Institute, Russian Academy of Sciences,
Gubkina str.~8, 119991 Moscow, RUSSIA; \\[1mm]
Department of Mechanics and Mathematics, Moscow Lomonosov State University,
Vorobiovy Gory, GSP-1, 119991 Moscow, RUSSIA

\medskip
\noindent
\emph{E-mail}:
\href{mailto:wadim@mi.ras.ru}{\texttt{wadim@mi.ras.ru}}

\end{document}